\magnification=\magstep1
\input amstex
\documentstyle{amsppt}
\NoRunningHeads
\TagsOnRight

\define\C{{\Bbb C}}
\define\R{{\Bbb R}}

\define\Dom{\operatorname{Dom}}

\define\E{{\Cal E}}
\define\X{{\Cal X}}
\define\Y{{\Cal Y}}

\define\Ker{\operatorname{Ker}}

\define\A{{\Cal A}}

\define\spec{\operatorname{spec}}
\define\sm{\operatorname{sm}}
\define\lar{\operatorname{la}}

\redefine\Im{\operatorname{Im}}

%\define\ln{\operatorname{ln}}
\redefine\H{\Cal H}

\define\La{\Lambda}
\define\la{\lambda}
\define\om{\omega}
\define\Om{\Omega}
\define\De{\Delta}
\define\de{\delta}
\define\Ga{\Gamma}
\define\ga{\gamma}
\define\U{{\Cal U}}

%\nopagenumbers

\topmatter
%\subheading{Preliminary version}

\title
DE RHAM THEOREM FOR EXTENDED $L^2$-COHOMOLOGY
\endtitle


\dedicatory
Dedicated to Selim Grigor'evich Krein  on the occasion of  his 80th birthday
\enddedicatory

\author
M.A.SHUBIN\footnotemark"$^1$"
\endauthor
\footnotetext"$^{1}$"{Partially supported
by US - Israel Binational Science Foundation Grant 94-00299/1}

\abstract
We prove an analogue of the de Rham  theorem for the extended $L^2$-cohomology
introduced by M.Farber
\cite{Fa}. This is done by establishing that the de Rham complex over a compact
closed manifold with
coefficients in a flat Hilbert bundle $E$ of $\A$-modules over
a finite von Neumann algebra $\A$
is chain-homotopy equivalent in the sense of [GS]
(i.e. with bounded morphisms and homotopy operators)
to a combinatorial complex with the same coefficients.
This is established by using the Witten deformation of the de Rham complex.
We also prove that the de Rham complex is chain-homotopy equivalent to the
spectrally truncated de Rham complex which is also finitely generated.
\endabstract

\affil
Northeastern University
\endaffil
\address
Department of Mathematics,
Northeastern University,
Boston, MA 02115
\endaddress
\email
shubin\@neu.edu
\endemail
%\smallarea{
%\centerline {\itsmall Department of Mathematics}
%\centerline {\itsmall Northeastern University,}
%\centerline {\itsmall Boston, MA 02115}
%\centerline {\itsmall e-mail shubin$\@$neu.edu}
%\centerline {\itsmall fax 617-373-5658}
%}
\endtopmatter

\document

\heading{Introduction}
\endheading
$L^2$-cohomology is a natural tool for constructing invariants of
non-compact manifolds.
To define it we need an extra  structure of the manifold near infinity, e.g.
a Riemannian metric or a triangulation, or rather a quasiisometry class of
metrics
or a class of triangulations which are compatible with a uniform structure.
Sometimes such an extra structure naturally exists, e.g. if the manifold is the
open set of all
regular points of an algebraic variety, or if it has a proper cocompact
action of
a discrete group (e.g. if it is a regular covering of a compact manifold).

Assume that we have such a manifold $Y$ and that an admissible
triangulation $T$ of it is given. Then we can define
Hilbert spaces of $L^2$-cochains $C^i_{(2)}=C^i_{(2)}(T,Y)$ and they form a
complex
$$\dots\longrightarrow C^{i-1}_{(2)}@>d_{i-1}>>C^i_{(2)}@>d_i>>C^{i+1}_{(2)}
\longrightarrow\dots\,,\tag 0.1$$
where $d_i$ are bounded linear operators.

There are two possible ways to define the $L^2$-cohomology of this complex.
One of them is just to form the usual cohomology
$$L^2H^i(T,Y)=\Ker d_i/\Im d_{i-1}\,,$$
ignoring topology. But then we have to face the fact that these
cohomology spaces are infinite-dimensional and
have no natural Hausdorff topology.

Another way is to form reduced cohomology
$$L^2H^i(T,Y)=\Ker d_i/\overline{\Im d_{i-1}}\,,$$
where the bar over $\Im d_{i-1}$ means its closure in the Hilbert space
$C^i_{(2)}$.
However, this leads to a big loss of information.

It was first noticed in \cite{NS}
that for covering manifolds  topological information can be extracted by
considering
the behavior
of the spectra of the Laplacians of the complex (0.1) near $0$. In
particular some
numbers characterizing the behavior of the spectrum near zero were defined
in \cite{NS1}
(they were later called Novikov-Shubin invariants). It was  proved in \cite{GS}
that they are in fact homotopy invariants. More details about these invariants
can be found in \cite{E, E1, GS, LL}. In particular it follows from
\cite{GS} that these
invariants are well-defined for arbitrary regular coverings of finite
CW-complexes.

M.Farber \cite{Fa} discovered a way to approach the spectrum-near-zero
phenomenon
through a new cohomology theory. It  is convenient to describe his approach
in a more
general context. Namely, let ${\tilde X}$ be a manifold with a free action of
a discrete group $\Ga$ so that the quotient manifold $X={\tilde X}/\Ga$ is
compact.
Then $L^2$-functions or $L^2$-forms on ${\tilde X}$ can be considered as
sections of an
infinite-dimensional vector bundle over $X$.
The fiber of this bundle is a finitely-generated
Hilbert module over the finite von Neumann algebra $\A$ associated with $\Ga$.

M.Farber actually considered even more general situation of an arbitrary
Hilbert bundle
over a finite CW-complex.  The fiber of this bundle should be
a finitely generated Hilbert $\A$-module, where $\A$ is an arbitrary finite
von Neumann
algebra.
In this setting M.Farber introduced an extended cohomology theory. It  takes
values in an abelian category $\E(\A)$ which is obtained by applying a
P.Freyd construction to
the additive category $\H(\A)$ of all finitely generated Hilbert $\A$-modules.
The category $\E(\A)$ has $\H(\A)$  as the full subcategory of all
projective objects.

Now a cochain complex in the category $\H(\A)$ can be considered as a
complex in $\E(\A)$
and the corresponding cohomology is called {\it extended  cohomology}. The
extended cohomology objects belong to the category  $\E(\A)$.

The  reduced $L^2$-cohomology spaces
appear then  as the projective parts of the extended cohomology objects.
But the extended cohomology objects contain
also the torsion parts, which determine  the Novikov-Shubin
invariants. There are also some other numerical
invariants determined by the torsion part of extended cohomology.
One of them is the number of generators, which was used in \cite{Fa}  to
improve the
von Neumann version of the Morse inequalities from \cite{NS}.

The construction suggested in \cite{Fa} used cochain complexes arising from
a cell decomposition of a  manifold or a more general polyhedron.
Our purpose in this paper is to define the de Rham version of the extended
$L^2$-cohomology of a manifold and to understand the corresponding de Rham
theorem.
This would simplify
calculations of  extended $L^2$-cohomology
and make possible its applications
in geometry and analysis.

In order to define the
de Rham version of  extended cohomology we will
consider spaces of $L^2$ differential  forms with values in a flat
vector bundle $E$ such that its fiber is a
 finitely generated Hilbert module over a
finite von Neumann algebra $\A$. These forms constitute a
twisted de Rham complex which is a complex of Hilbert
$\A$-modules, and the differentials are
(generally unbounded) closed, densely defined linear operators. Unfortunately,
the Hilbert $\A$-modules of forms are not finitely generated, and it does not seem
possible to generalize the construction of the extended category \cite{Fa} to
include such complexes. In order to overcome this difficulty, we
consider  the homotopy type of the twisted de Rham complex with respect
to chain-homotopy equivalence as studied in \cite{GS, GS1} (i.e. with bounded
morphisms and homotopy operators). Then we define {\it extended de Rham
cohomology}
as the extended cohomology of any finitely generated Hilbert $\A$-complex
which is
homotopy-equivalent to the de Rham complex. Such a finitely generated
$\A$-complex is called a {\it finite approximation} of the de Rham complex.
The result does not depend on the choice of  approximation.

Now we have to connect the extended de Rham cohomology with  extended
combinatorial chain cohomology. We prove that they coincide and this is what
we call the {\it de Rham theorem in extended cohomology}. In fact the de Rham
complex is homotopy-equivalent to the corresponding combinatorial complex,
so the combinatorial cochain complex is a finite approximation
of the de Rham complex. This is proved with the use of the Witten deformation
of the de Rham complex. It is important here that  in the Witten
approximation we have
 a gap
separating the ``small eigenvalues" from the rest of the spectrum (see
\cite{BFKM} and \cite {Sh2}).

Another finite approximation of the de Rham complex can be obtained
by truncating all the $L^2$-form
spaces by the spectral projections of the Laplacians. So the desired
approximation is just the subcomplex corresponding to the {\it ``small
eigenvalues"}.
It is more difficult here to establish that this ``small eigenvalues" complex
is really finitely  generated. To do this we use uniform  kernel estimates and
factor decomposition of the fiber.

The de Rham type theorems which are established in this paper generalize
two well-known
theorems of this kind: a result of  J.Dodziuk \cite{D}, who
studied the  reduced  $L^2$-cohomology, and a result of
A.Efremov \cite{E, E1} who proved the equality of the combinatorial and
analytic
Novikov-Shubin invariants.

W.L\"uck \cite{L\"u} suggested a more algebraic cohomological approach to the
description of the
spectra near zero which is equivalent to Farber's approach.  L\"uck's technique
can be used in the De Rham theory as well, though Farber's approach seems more
convenient.

After becoming acquainted with my idea of  de Rham theory in extended
cohomology
M.Farber suggested a different approach to de Rham theory. His approach leads
to a natural relation of  extended de Rham theory with  \v Cech cohomology
and makes use of sheaf theory (\cite{Fa1}).

I am grateful to M.Farber, V.Mathai, R.Nest and A.Suciu for useful discussions.

\heading{1. Preliminaries}
\endheading
In this section we will
briefly describe Farber's  extended cohomology theory \cite{Fa}.

Let $\A$ be a finite von Neumann algebra. This means that $\A$ is a von Neumann
algebra which has a finite  trace $\tau:\A\to\C$ which is faithful
and normal. We will fix this trace and we will
always assume it normalized, i.e. $\tau(1)=1$.
If $\A$ happens to be a factor, then it has to be of type $\hbox{I}_n$
 or $\hbox{II}_1$. More generally a finite von Neumann algebra $\A$ is a
direct integral of type $\hbox{I}_n$
 and $\hbox{II}_1$ factors.
A  definition of  finite von Neumann algebra in terms of projections
can be found in \cite{T} (Definition 1.16, Ch.V, p.296).

Recall that $\A$ can be
 represented
as a weakly closed selfadjoint subalgebra of the algebra of all bounded
linear operators in a Hilbert space. The operator norm   induces a norm
in $\A$ which does not depend on the choice of the representation. We will refer
to this norm as the $C^\ast$-norm in $\A$. The corresponding topology in $\A$
will be called the {\it norm topology}.

The operation $B\to B^\ast$ (taking the adjoint operator)  is a well-defined
involution in $\A$.

Let $\ell^2(\A)$ denote the Hilbert space obtained as the
completion of $\A$ with respect to the inner product given by the
trace: $(a,b)=\tau(b^\ast a)$ for all $a,b\in \A$.

A {\it Hilbert module} over $\A$
is a Hilbert space $M$
together with a continuous left $\A$-module structure
(here $\A$ is considered with its norm topology) such that there
exists an isometric $\A$-linear embedding of $M$ into
$\ell^{2}(\A)\otimes H$, for some Hilbert space $H$.
A Hilbert module $M$ is {\it finitely generated} if it
admits an embedding $M\rightarrow\ell^{2}({\Cal A})\otimes H$ as above with
finite-dimensional $H$.

It is often more convenient to forget about the inner product in a Hilbert
module $M$ and consider this module up to an $\A$-linear topological
isomorphism.
In \cite{Fa} the corresponding class is called
 {\it  Hilbertian module}. In other words a Hilbertian module
is a topological vector space $M$ with a structure of a left (algebraic)
$\A$-module such that
the action of $\A$ is continuous and there exists an inner product
$(\;,\;)$ on $M$ which generates the topology of $M$ and such
that $M$ together
with $(\;,\;)$ and with the $\A$-action is a Hilbert
module. The difference between Hilbert and Hilbertian modules is important
in categorical
constructions but it is not important for us. Therefore we shall only deal
with Hilbert
$A$-modules, though isomorphisms of Hilbert $\A$-modules will be understood
as arbitrary
$\A$-linear topological isomorphisms (not necessarily unitary ones).

Let $\H(\A)$ denote the additive category whose objects
are finitely generated left Hilbert $\A$-modules  and whose morphisms are
continuous $\A$-module  homomorphisms.
This category  depends on the choice of the trace
$\tau$.

The P.Freyd construction \cite{Fr} provides an extended category $\E(\A)$
which is an abelian category, containing  $\H(\A)$ as a full subcategory.
An {\it object} of the category $\E(\A)$
is defined as a
morphism $(\alpha:A^\prime\to A)$ in the category $\H(\A)$.
Given a pair of objects $\X=(\alpha: A^\prime \to A)$ and
$\Y=(\beta:B^\prime\to B)$ of $\E(\A)$,
a {\it morphism} $\X\to\Y$ in the category $\E(\A)$ is an
equivalence class of morphisms $f:A\to B$ of the category $\H(\A)$
such that $f\circ\alpha=\beta\circ g$
for some morphism $g:A^\prime \to B^\prime$ in $\H(\A)$.
Two morphisms $f:A\to B$ and $f^\prime:A\to B$ of $\H(\A)$ represent
the same morphism $\X\to\Y$ in $\E(\A)$
 iff
$f-f^\prime = \beta\circ F$ for some morphism $F:A\to B^\prime$ of the category
$\H(\A)$. This defines an equivalence relation.
%The morphism $\X\to\Y$, represented by $f:A\to B$,
%is denoted
%$$[f]:(\alpha:A^\prime\to A)\ \to\ (\beta:B^\prime\to B)\quad
%\text{or by}\quad
%[f]:\X\to\Y.$$
The {\it composition} of morphisms is defined as the composition of the
corresponding morphisms $f$ in the category $\H(\A)$.

Since $\E(\A)$ is an abelian category,
any cochain complex $D$ in $\H(\A)$ (or in $\E(\A)$)
has cohomology $\H^\ast(D)$ which is
well-defined as  a graded object of $\E(\A)$. It is called the {\it
extended cohomology}
and it has all the
standard cohomological properties.

In particular it is functorial, i.e.
any chain map between cochain complexes $g: D_1\to D_2$ in $\H(\A)$
induces a morphism of their extended cohomologies
\newline
$g^\ast: \H^\ast(D_1) \to \H^\ast(D_2)$.

Extended cohomology  is
homotopy invariant, i.e. the induced morphism $g^\ast$ depends only
on the homotopy class of $g$.

There is a long exact sequence
$$\dots \to \H^i(D_1) \to \H^i(D_2)\to \H^i(D_3) \to
\H^{i+1}(D_1)\to \dots$$
corresponding to any short exact sequence
$$0\to D_1\to D_2\to D_3\to 0$$
of cochain complexes in $\H(\A)$.

\medskip
Now let us recall the definition of extended cohomology in a topological
context.

Let $\Ga$ be a countable discrete group and
$M$ be a finitely generated Hilbert $\A$-module. Suppose that we have a
right action of $\Ga$
on $M$ by $\A$-automorphisms i.e.  a representation
$\rho:\Ga^{op}\to GL_{\A}(M)$. Here
$GL_{\A}(M)$ is the group of all (bounded) $\A$-automorphisms of $M$
(not necessarily unitary), and
$\Ga^{op}$ is the group which is obtained from $\Ga$
by introducing a new operation so that the product of $\alpha$ and $\beta$
in $\Ga^{op}$ equals to the product $\beta\alpha$ in $\Ga$. Following
\cite{BFKM}
we will say in this case  that $M$ is a finitely generated Hilbert
$(\A,\Ga^{op})$-module.

Note that by definition
the actions of $\A$ and $\Ga^{op}$ on $M$ commute.

The simplest important example of this situation appears when
$M=\ell^2(\Ga)$, and
$\A={\Cal N}(\Ga)$ is the von Neumann algebra generated by the left
translations by
the elements $\ga\in\Ga$ in $\ell^2(\Ga)$. There is also a natural right
action of $\Ga^{op}$
in $\ell^2(\Ga)$ by right translations. This action obviously commutes with
the action of
${\Cal N}(\Ga)$. Therefore $\ell^2(\Ga)$ has the structure of a finitely
generated Hilbert
$({\Cal N}(\Ga),\Ga^{op})$-module.

Let $X$ be a finite, connected CW-complex and $\Ga$ its fundamental group.
Let $\tilde X$ be  the universal covering of $X$.
The group $\Ga$ acts on $\tilde X$ by deck transformations. Therefore
for any finitely generated Hilbert $(\A,\Ga^{op})$-module $M$ we can form
the space $E=M\times_\Ga \tilde X$. There is a natural projection $p: E\to X$
so that $E$ becomes a flat Hilbert $\A$-bundle with  standard fiber $M$
and representation $\rho$ giving the monodromy.

A  flat Hilbert $\A$-bundle $E$ with  fiber $M$ over a topological
space $X$ can also be defined by a 1-cocycle on $X$ with values in
$GL_{\A}(M)$.
Namely, for a sufficiently fine open covering $\U$ of $X$ and any
$(U,V)\in\U\times\U$ we should be given an element $g_{UV}\in GL_{\A}(M)$
so that the
following conditions are satisfied:

 (a) $g_{UU}=1$;

(b) $g_{UV}\cdot g_{VW}=g_{UW}$ if $U\cap V\cap W\not=\emptyset$.

The total space of the bundle
$E$ is constructed then by identifying the points
$(x,m)\in U\times M$ with $(x, g_{VU}\cdot m)\in V\times M$
in the disjoint union of the spaces $U\times M$ for $x\in U\cup V$.

Let $X^\prime\subset X$. Then  any flat Hilbertian bundle $E$
over $X$ has a naturally defined restriction
$E^\prime=E|_{X^\prime}$ which is a flat Hilbertian
bundle over $X^\prime$.

Assume for a moment that $X$ is a compact closed $C^\infty$-manifold.
Denote by $C^\infty(X,E)$ and $L^2(X,E)$ the spaces of $C^\infty$-sections and
$L^2$-sections of $E$, respectively. It is easy to see that if we take
$M=\ell^2(\Ga)$ with the action of $\Ga^{op}$ by right translations as described
above, then the space $L^2(X,\ell^2(\Ga))$ is naturally identified with the
space $L^2(\tilde X)$ defined by a $\Ga$-invariant measure (i.e. a measure
lifted from $X$) with a smooth positive density. Under this
identification the space $C^\infty(X,\ell^2(\Ga))$ becomes the space
$H^\infty(\tilde X)$ which can be described as the intersection of all uniform
Sobolev spaces on $\tilde X$ (see \cite{Sh1} for more details about these
spaces)
or as the space of all $u\in C^\infty(\tilde X)$ such that
$Du\in L^2(\tilde X)$ for every $\Ga$-invariant differential operator $D$
with $C^\infty$
coefficients.

Let us return to the general case when $X$ is a finite, connected CW-complex.
Then the  universal covering $\tilde X$ has the natural structure
of a CW-complex generated by the lifts of the cells of $X$.
A flat bundle $E$ as above defines an associated
locally constant sheaf $\underline {E}$. The space of  sections of
$\underline E$
over a given cell is isomorphic to the standard fiber $M$. Using the
restriction maps
in the usual way,
we obtain a finitely generated Hilbert cochain $\A$-complex which we will denote
by $({\Cal C}^\bullet,\de)$. It can be also defined purely algebraically
(see \cite{Fa}
for the definitions in  this particular context).

The extended cohomology of $({\Cal C}^\bullet,\de)$ does not depend on the
choice of  CW representation of $X$ (\cite{Fa}). We will call it the
{\it extended cohomology of} $X$ (with coefficients in $E$) and
denote it by  ${\Cal H}^\ast(X,E)$.

In the case  $X$ is a compact closed manifold, we will actually use
a special CW-complex which is associated with a generic Morse function
as described in  F.Laudenbach's appendix to \cite{BZ} and also in \cite{BFKM}.

\heading{2. Fredholm complexes and their extended $L^2$-cohomology}
\endheading
We would like to generalize  the definition of  extended cohomology to
a category of  Fredholm complexes of Hilbert
$\A$-modules with unbounded differentials (like the one  considered in
\cite{GS}).
The extended cohomology objects will still belong to  the extended category
$\E(\A)$.

Let
$$C:\ \ \dots \to C^{i-1}@>d_{i-1}>>C^i@>d_i>>C^{i+1}\to \dots$$
be a sequence of finite length formed by Hilbert $\A$-modules and closed,
densely defined linear operators $d_i$. It is called a Hilbert cochain
$\A$-complex if
the following two conditions are satisfied:

(i) $d_{i+1}d_i=0$ on $\Dom(d_i)$ for all $i$
(here $\Dom(d_i)$ denotes the domain of $d_i$);

(ii)  For every $i$ the operator $d_i$ commutes with  multiplication by any
$a\in\A$
in the sense that $d_i a=ad_i$ on $\Dom(d_i)$.

In particular we require that $a(\Dom(d_i))\subset\Dom(d_i)$ for all
$a\in\A$ and all $i$.

Assume that we have two Hilbert cochain $\A$-complexes $C$ and $C_1$ with
the differentials $d_i$, $d^\prime_i$ respectively. A {\it morphism of Hilbert
cochain $\A$-complexes} $f:C\to C_1$ is  a set of {\it bounded} linear
$\A$-morphisms
$f_i:C^i\to C_1^i$ such that $d^\prime_if_i=f_{i+1}d_i$ on $\Dom(d_i)$.
With this set of morphisms the Hilbert cochain $\A$-complexes form a category.

Two morphisms $f,g:C\to C_1$ of Hilbert $\A$-complexes are called {\it
homotopic}
if there exist bounded homotopy
operators ($\A$-morphisms) $T_i:C^i\to C_1^i$ such that
$f_i-g_i-d_{i-1}T_i=T_{i+1}d_i$ on $\Dom(d_i)$ for all $i$. The homotopy is
an equivalence
relation on the set of all morphisms. Hence the notion of {\it
chain-homotopy equivalence}
of Hilbert cochain $\A$-complexes is well defined.

\medskip
{\bf Definition 2.1.} We will say that a Hilbert $\A$-complex $C$ is {\it
Fredholm}
if there exists a Hilbert cochain $\A$-complex $D$
of finite length such that the following two conditions are satisfied:

(i) all $D_i$ are finitely generated Hilbert $\A$-modules and all the
differentials are bounded
(i.e. $D$ is a complex in the category $\H(\A)$);

(ii) $D$ is chain-homotopy equivalent to $C$ in the sense described above.

 A pair $(D, f)$, consisting of a chain complex
$D$ of finite length in $\H(\A)$ and of a chain-homotopy equivalence
$f:D\to C$ will be called a {\it finite approximation of $C$}.

For a Fredholm Hilbert $\A$-complex $C$ we may define its {\it extended
cohomology} $\H^\ast(C)$
as the extended cohomology $\H^\ast(D)$, where $(D, f)$ is an arbitrary
finite approximation of $C$.

Suppose that we have a cochain map $g:C_1\to C_2$ between two
Fredholm Hilbert cochain $\A$-complexes.  Let
$$f_1:D_1\to C_1\quad\text{and}\quad f_2:D_2\to C_2$$
be  finite approximations. Then the composition
$$D_1@>{f_1}>>C_1@>g>>C_2@>h>>D_2,$$
where $h$ is a homotopy inverse to $f_2$, is a chain map between two
cochain complexes in $\H(\A)$. Thus there is well-defined induced map
$\H^i(D_1)\to \H^i(D_2)$ which we understand as the morphism
induced by the original cochain map $g$; sometimes we will denote this map
$g^\ast:\H^i(C_1)\to \H^i(C_2)$.

In particular we see that there exists a canonical isomorphism
$\H^i(D_1)\to \H^i(D_2)$ if $D_1$ and $D_2$ are two
arbitrary finite approximations of a given Fredholm complex $C$. Thus,
the extended cohomology of a finite approximation of a Fredholm complex
does not depend on the approximation.

In the geometric situation, the complex $C$ will arise as a
de Rham complex of $L^2$-differential forms on a compact smooth
manifold (without boundary)
with values in a flat Hilbert bundle over $\A$. We will give constructions of
two finite approximations of  this de Rham complex
(constructions of this kind we will call de Rham
theorems). In particular, we will establish that the
finite cochain complex in $\H(\A)$ constructed by using a cell decomposition
of the
manifold (which was used in \cite{Fa} to define the extended cohomology)
is chain-homotopy equivalent to the de Rham complex; this will be our
version of the de Rham
theorem for  extended $L^2$-cohomology.

\heading{3. De Rham complex}
\endheading
Let $\A$ be a finite von Neumann algebra with a fixed trace $\tau$,
and $E$ a fixed Hilbert $\A$-bundle over a compact closed  manifold $X$.
Let $M$ be the standard fiber of $E$, so $M$ is
a fixed finitely generated Hilbert $\A$-module.

Let $\{g_{UV}\}$ be a 1-cocycle (with values in $GL_{\A}(M)$) which defines $E$.

For any open set $U\subset X$
denote by $\La^p(U)$ the set of all smooth $p$-forms on $U$, and
$\La^\bullet(U)=\oplus_p\La^p(U)$.

Denote by $C^\infty(U,M)$ the  Fr\'echet space of smooth functions on $U$
with values in $M$. The smooth $p$-forms on $U$ with values in $M$ form a
Fr\'echet space
$$\La^p(U,M)=\La^p(U)\otimes_{C^\infty(U)}C^\infty(U,M)\,.$$
The Fr\'echet space of smooth $p$-forms with values in $E$ can be defined as
$$\La^p(X,E)=\Lambda^p(X)\otimes_{C^\infty(X)}C^\infty(X,E)\;,$$
where $C^\infty(X,E)$ is the space of $C^\infty$-sections of $E$ over $X$.
In other words a form $\om\in\La^p(X,E)$ is a collection
$\om=\{\om_U\}_{U\in\U}$
where $\om_U\in\La^p(U,M)$, such that for any pair $U,V\in\U$ with $U\cap
V\not=\emptyset$
$$\om_U|_{U\cap V}\ =\ g_{UV}\cdot \om_V|_{U\cap V}.\tag 3.1$$
The exterior differential naturally extends to the operator
$$d=d^{E}:\La^p(X,E)\longrightarrow \La^{p+1}(X,E),\tag 3.2$$
which is defined by applying the usual de Rham differential $d$ to all $\om_U$
(this is well-defined because the functions $g_{UV}$  in the cocycle
condition (3.1) are
constant).

 The operator $d^{E}$ is also called the {\it covariant derivative} in $E$.
Clearly
$(d^{E})^2=0$, so we have a well-defined de Rham complex
$$0@> >>\La^0(X,E)@>d>>\dots@>d>>\La^p(X,E)
@>d>>\La^{p+1}(X,E)@>d>>\dots@>d>>\La^n(X,E)
@> >> 0\;,\tag 3.3$$
where $n=\dim_{\R}X$.

Let us choose a Riemannian metric on $X$ and a smooth Hermitian metric on
the bundle $E$,
such that the topology it induces in each fiber coincides with the standard
topology
of that fiber. Then we can define the space $L^2\La^p(X,E)$ of $E$-valued
$L^2$-forms of
degree  $p$. This space does not depend on the choices of the Riemannian
metric on $X$
and the Hermitian metric on $E$. Now we can define the $L^2$-de Rham complex
$$0@> >>\dots@> >>L^2\La^p(X,E)
@>d>>L^2\La^{p+1}(X,E)@> >>\dots@> >> 0\;,\tag 3.4$$
where the differentials are closures of the corresponding differentials in
(3.3).
Clear--ly this complex is a Hilbert $\A$-complex in the sense described in
Sect.2.

Now we can formulate our de Rham theorem. A similar theorem for smooth forms
has been proved by M.Farber \cite{Fa1}.

\proclaim{Theorem 3.1} Let $X$ be a compact closed $C^\infty$-manifold
without boundary.
Then the de Rham complex (3.4) is Fredholm in the sense of Definition 2.1.
Its extended cohomology coincides with the combinatorially defined extended
cohomology
of  $X$ with coefficients in $E$.
\endproclaim

In the next two sections we will give two constructions of  finite
approximations for
the complex (3.4). In particular we will prove that this complex is
homotopy equivalent
to the combinatorial cochain $L^2$-complex.

\heading {4. Witten approximation}
\endheading

 We shall use the same notations
as in the previous section. In particular $X$ will be a compact closed
$C^\infty$-manifold.

Let $f:X\to\R$ be a Morse function on $X$. Consider the Witten deformations
of the de Rham
complexes (3.3) and (3.4) which are obtained by replacing $d$ by the
deformed differential
$$d_t=e^{-tf}de^{tf}=d+tdf\wedge\cdot,\ \ t\in\R\;.$$
In this way we obtain new complexes
$$0@> >>\dots@> >>\La^p(X,E)
@>d_t>>\La^{p+1}(X,E)@> >>\dots@> >> 0\;,\tag 4.1$$
and
$$0@> >>\dots@> >>L^2\La^p(X,E)
@>d_t>>L^2\La^{p+1}(X,E)@> >>\dots@> >> 0\;.\tag 4.2$$

Assume now that a Riemannian metric $g$ is chosen on $X$.
Then we can also define the corresponding Witten Laplacian
$$\De_t=d_t^\ast d_t+d_td^\ast_t\,.\tag 4.3$$
It acts in $\La^p(X,E)$ for every $p$ and also defines a self-adjoint
operator in
$L^2\La^p(X,E)$.

Following [BFKM]  we will assume for simplicity that the pair $(f,g)$ forms a
 {\it generalized triangulation}. This means that the following conditions
are satisfied:

(i) $f$ is self-indexing i.e. $f(x)=\hbox{index}(x)$ for any critical point
$x$ of $f$;

(ii) in a neighborhood of any critical point $x$ of $f$ there exist coordinates
$y_1,\dots, y_n$ such that in these coordinates $f$ has the form
$$f(y)=k-(y_1^2+\cdots +y_k^2)/2+(y_{k+1}^2+\cdots +y_n^2)/2\,,$$
where $k=\hbox{index}(x)$;

(iii) for any two critical points $x$ and $y$ of $f$, the unstable manifold
$W_x^-$
and the stable manifold $W_y^+$, associated to the vector field
$\hbox{grad}_g f$ intersect
transversally.

\proclaim{ Lemma 4.1} {\rm ([BFKM])} There exist positive constants
$C^\prime, C^{\prime\prime}$ and $t_0$ such that for any $t\geq t_0$
$$\spec(\De_t)\cap (e^{-tC^\prime}, C^{\prime\prime}t)=\emptyset\,,$$
where $\operatorname{spec}$ means the spectrum in $L^2\La^\bullet(X,E)$.
\endproclaim

{\bf Remark.} A stronger statement about  the limit behavior of the whole
spectrum
was independently proved in [Sh2]. In [Sh2] only a particular case of the
von Neumann
algebra and the Hilbert module associated with a regular covering of $X$
is considered but the arguments are easily extended to the general situation.

Lemma 4.1.allows to define  a splitting of the deformed de Rham complex
(4.2) into
a direct sum of two Hilbert $\A$-subcomplexes:
$$L^2\La^\bullet(X,E)=L^2\La^\bullet_{\sm}(X,E)\oplus
L^2\La^\bullet_{\roman la}(X,E),\ \ t>t_0\,.\tag 4.4$$
Here $L^2\La^\bullet_{\sm}(X,E)$ and $L^2\La^\bullet_{\lar}(X,E)$
are the images of the spectral projections of $\De_t$ corresponding to
``small" and ``large" eigenvalues (i.e. to the spectral intervals $[0,1)$
and $[1,\infty)$,
respectively). Denote these projections by $P_{\sm}(t)$ and $P_{\lar}(t)$.
They commute with the action of $\A$. They also commute  with $d_t$ because
$\De_t$ commutes with $d_t$. Therefore both  projections $P_{\sm}(t)$ and
$P_{\lar}(t)$ are morphisms of Hilbert cochain complexes.

Of course the decomposition
(4.4)  depends on $t$ but we do not indicate this explicitly so as to simplify
the notations.

All the complexes in (4.4) are considered with the differential $d_t$.

Denote by
$i_{\sm}(t): L^2\La^\bullet_{\sm}(X,E)\to  L^2\La^\bullet(X,E)$
the canonical inclusion. It is also a morphism of Hilbert cochain
$\A$-complexes.

\proclaim{Theorem 4.2} {\rm (i)} The maps $i_{\sm}(t)$ and $P_{\sm}(t)$
define a chain-homotopy equivalence of the Hilbert cochain $\A$-complexes
$L^2\La^\bullet_{\sm}(X,E)$ and $ L^2\La^\bullet(X,E)$ (with the
differentials $d_t$) for any $t>t_0$.

{\rm (ii)} The differentials $d_t$ are bounded $\A$-operators in
 $L^2\La^\bullet_{\sm}(X,E)$.

{\rm (iii)}
 $L^2\La^\bullet_{\sm}(X,E)$  is a finitely generated
Hilbert cochain $\A$-complex.
\endproclaim

\proclaim{Corollary 4.3} The de Rham complex  $L^2\La^\bullet_{\sm}(X,E)$
is a Fredholm Hilbert cochain $\A$-complex in the sense of Definition 2.1.
\endproclaim

So Theorem 4.2. gives us the first  finite approximation of the de Rham
complex.

\medskip\noindent
{\bf Proof of Theorem 4.2.} We have
$$\|d_t\om\|^2=(d_t^\ast d_t\om,\om)\leq(\De_t\om,\om)\leq\|\om\|^2,
\ \om\in L^2\La^\bullet_{\sm}(X,E)\,,$$
which proves (ii).

To prove (i) let us introduce a Green operator $G_t$ by the formula
$$G_t=\De_t^{-1}(I-P_{\sm}(t)): L^2\La^\bullet(X,E)
\to L^2\La^\bullet(X,E)\,,$$
where $I$ is the identity operator, $\De_t^{-1}$ is the operator
which is inverse to $\De_t$ restricted to
$(L^2\La^\bullet_{\sm}(X,E))^\perp=\Im(I-P_{\sm}(t))$, and
by definition $G_t=0$ on $L^2\La^\bullet_{\sm}(X,E)=\Im P_{\sm}(t)$.
Clearly $G_t$ is a bounded $\A$-operator. It follows easily from  elliptic
regularity
arguments that $G_t$ is smoothing by 2 units in the corresponding Sobolev scale
(see e.g. similar arguments in [Sh]). In particular, the operators
$d_tG_t$, $d_t^\ast G_t$,
$d_t^\ast d_tG_t$ and $d_td_t^\ast G_t$ are  bounded in $L^2\La^\bullet(X,E)$.

Note also that $G_t$ commutes with $d_t$ and $d_t^\ast$.  Therefore we obtain
$$I-P_{\sm}(t)=\De_t\De_t^{-1}(I-P_{\sm}(t))=
d_td_t^\ast G_t+d_t^\ast d_tG_t=d_t(d_t^\ast G_t)+(d_t^\ast G_t)d_t\,,$$
so the operator $d^\ast_t G_t$ supplies a chain homotopy between the
$\A$-endomorphisms
$I$ and $P_{\sm}(t)$ of the Hilbert cochain $\A$-complex $L^2\La^\bullet(X,E)$.
This proves (i).

To prove (iii) we use Theorem 5.5 from [BFKM] which provides an isomorphism
of Hilbert
cochain $\A$-complexes
$$\phi:(L^2\La^\bullet(X,E),\tilde d_t)\longrightarrow ({\Cal C}^\bullet, \de)
\,,\tag 4.5$$
where $\tilde d_t$ is obtained from $d_t$ by scaling
$$\tilde d_t=e^t\Bigl({t\over \pi}\Bigr)^{-1/2}d_t\,,$$
and $({\Cal C}^\bullet, \de)$ is a finite combinatorial $\A$-complex
associated with
the given bundle $E$ and the chosen generalized triangulation (see Section
4 in [BFKM]
and also \cite{BZ1} for the case $\A=\C$).
Since the complex  $({\Cal C}^\bullet, \de)$ is obviously finitely
generated, we obtain (iii).
$\square$

\proclaim {Corollary 4.4}
The extended de Rham cohomology with coefficients in $E$ coincides with
extended
combinatorial cohomology with the same coefficients.
\endproclaim

{\bf Proof.} The $L^2$-de Rham complex $(L^2\La^\bullet(X,E), d)$ is
obviously isomorphic
to its Witten deformation $(L^2\La^\bullet(X,E), d_t)$, so they have the same
extended cohomology. On the other hand,  replacing $d_t$ by $\tilde d_t$ also
leads to an isomorphic complex. It follows that the combinatorial complex
$({\Cal C}^\bullet, \de)$ is a  finite approximation for the de Rham
complex, therefore
it can be used to calculate its extended cohomology.
$\square$

{\bf Proof of Theorem 3.1.} According to  Theorem 6.4 in [Fa] we can use any
CW-complex structure on $X$ to calculate the  extended cohomology. In particular
we can use the above-mentioned complex $({\Cal C}^\bullet, \de)$. $\square$

{\bf Remark.} We see that both $({\Cal C}^\bullet, \de)$ and
$(L^2\La_{\roman sm}^\bullet(X,E), d_t)$ provide finite approximations of
the de Rham
complex $(L^2\La^\bullet(X,E), d)$. Another finite approximation will be
given in the
next section.

\heading {5. Spectral approximation}
\endheading
In this section we will show that it is possible to construct a finite
approximation
of the de Rham $\A$-complex (3.4) by the direct spectral cut-off i.e. without
first deforming  as above.

%Here is our second de Rham type theorem.

As before let $X$ be a closed smooth manifold and $E$
 a flat Hilbert $\A$-bundle over $X$.
Choose a Riemannian metric on $X$ and a Hermitian metric on $E$. Then the
Hilbert
spaces of $E$-valued forms $L^2\La^p(X,E)$ and $L^2\La^\bullet(X,E)$ are well
defined. We can also consider  the Laplacian
$$\De^{(p)}=d^\ast d+dd^\ast:\La^p(X,E)\longrightarrow \La^p(X,E)$$
which is a self-adjoint
operator in  $L^2\La^p(X,E)$. Denote also $\De=\oplus_p\De^{(p)}$, so $\De$
acts in
$\La^\bullet(X,E)$ and is a
self-adjoint operator in $L^2\La^\bullet(X,E)$.

Denote by $E_\lambda^{(p)}$  the spectral projection of the Laplacian
$\Delta^{(p)}$ corresponding to the spectral interval $[0,\lambda]$ where
$\lambda>0$. Denote
$L^{(p)}_\la=\Im E_\lambda^{(p)}$, so  $L^{(p)}_\la$ is a closed subspace
in $L^2\La^p(X,E)$.

%\newpage
\proclaim{Theorem 5.1}

{\rm (i)} $L^{(p)}_\la\subset \La^p(X,E)$ (i.e. $L^{(p)}_\la$ consists of
smooth $E$-valued forms) and $L^{(p)}_\la$ is a Hilbert  $\A$-submodule in
$L^2\La^p(X,E)$\,.

{\rm (ii)}
$d(L^{(p)}_\la)\subset L^{(p+1)}_\la$\; and $d:L^{(p)}_\la\to L^{(p+1)}_\la$ is
a bounded $\A$-operator.

{\rm (iii)} The inclusion
$i_\lambda:(L^{(\bullet)}_\la,d)\to (L^{(\bullet)},d)$
is a  chain-homotopy equivalence in the sense of section 2; the  homotopy
inverse map
is given by the projection $E_\la^{(p)}$.

{\rm (iv)} The Hilbertian $\A$-module  $L^{(p)}_\la$  is finitely generated
for every $p$ and every $\lambda>0$\;.
\endproclaim
\newpage

Thus the complex $(L^{(\bullet)}_\la,d)$ corresponding to ``small eigenvalues"
 provides another finite approximation of the de Rham complex.

\subheading{ Proof of (i)-(iii) in Theorem 5.1} The proof is done by the same
arguments as  the proof of (i) and (ii) in Theorem 4.2. Namely, (i) follows
from  elliptic regularity, the proof of (ii) is identical to the proof of (ii)
in Theorem 4.2, and the proof of (iii) uses a homotopy operator constructed
from a  Green operator as in the proof of (i) in Theorem 4.2.
$\square$
\medskip
To prove (iv) we will first consider the simplest case.

\subheading {Proof of (iv) when $\A$ is a factor}
We have to prove that $L^{(p)}_\la$ is finitely generated. We shall
use the natural extension of the trace $\tau$ to morphisms of
Hilbert $\A$-modules. By an abuse of notation we will also denote this
extension by $\tau$.
Note first that $\tau(E_\lambda^{(p)})<\infty$ due to arguments similar to
the ones given in \cite{A}.

The assumption that $\A$ is a factor allows us to use the following
simplifying fact:
 for any two Hilbert $\A$-modules $L_1, L_2$ a continuous inclusion
of $\A$-modules $L_1\subset L_2$ exists if and only if
$\dim_\tau L_1\leq \dim_\tau L_2$ where $\dim_\tau$ is the von Neumann
dimension function induced by the trace $\tau$. Since
$\dim_\tau L^{(p)}_\la=\tau(E_\lambda^{(p)})<\infty$, the Hilbert
$\A$-module $L^{(p)}_\la$ can be embedded as a Hilbert $\A$-submodule
into a free Hilbert $\A$-module $l^2(\A)\otimes\C^N$,
 for a sufficiently large integer $N$. This  means by definition that
$L^{(p)}_\la$ is finitely generated. $\square$

The general case requires
 some preparations concerning direct integral factor decompositions
of the algebra $\A$ and finitely generated Hilbert $\A$-modules.

In the general case $\A$ is a direct integral of finite factors.
 We shall deduce the proof in this case  from the case when
$\A$ is a factor.
Let us start with a description of the finitely generated Hilbert
$\A$-modules in terms of their direct integral factor decomposition.

We refer to \cite{Di,P,Sc}  for  results about the direct integral
decompositions of von Neumann algebras and traces.

Consider the direct integral factor decomposition of $\A$:
$$\A=\int_\Om\oplus\A(\om)\mu(d\om)\tag 5.1$$
Here $\Om$ is a Borel space, $\om\mapsto\A(\om)$ is a Borel
family of finite factors, $\mu(d\omega)$ is a finite Borel measure on
$\Omega$. We will choose the Borel space $\Om$ and the measure $\mu(d\om)$
in a special way: we will assume that $\A$ has a direct integral decomposition
(5.1) in its regular representation, i.e. in its natural representation in
$\ell^2(\A)$.
So we have
$$\ell^2(\A)=\int_\Om\oplus h(\om)\mu(d\om)\,,\tag 5.2$$
where $h(\om)$ is a Hilbert space where $\A(\om)$ is represented.

Now let us consider an arbitrary finitely generated Hilbert $\A$-module $M$.

\proclaim{Lemma 5.2} $M$ has  a direct
integral decomposition
$$M=\int_\Om M(\om)\mu(d\om), \tag 5.3$$
where  $M(\om)$ is a Hilbert $\A(\om)$-module and in this decomposition the
action of $\A$ is diagonal  in the obvious sense with respect to the
decomposition (5.1).
\endproclaim
\newpage
{\bf Proof.} Let  $M$ be a Hilbert submodule in
$(\ell^2(\A))^N=\ell^2(\A)\otimes\C^N$, and $P$
be the orthogonal projection in $(\ell^2(\A))^N$ with  $\Im P=M$.
Clearly $P$ should be
in the commutant of the diagonal action of $\A$ in $(\ell^2(\A))^N$. According
to Theorem 4.11.8 in [P] this commutant is again the direct integral
$$\A^\prime=\int_\Omega\oplus\A^\prime(\omega)\mu(d\omega)\tag 5.4$$
with the diagonal action in $(\ell^2(\A))^N$. Therefore we have
$$P=\int_\Om\oplus P(\om)\mu(d\om)\,,\tag 5.5$$
where $P(\om)$ is an orthogonal projection in $h(\om)$ such that
$P(\om)\in \A^\prime(\om)$.
Denote $M(\om)=\Im P(\om)\subset h(\om)\otimes\C^N$. Then (5.5) implies
the decomposition (5.2). $\square$

\medskip
The trace $\tau$ on $\A$ can be written in the form
$$\tau(a)=\int_\Omega \tau_\omega(a(\omega))\mu(d\omega),
\ \ a\in\A,
\tag 5.6$$
where
$a=(a(\omega))$ in the direct integral decomposition (5.1), and
$\tau_\omega$ is a trace on $\A(\omega)$.

Note that we do not assume the traces $\tau_\omega$ to be normalized. For
future use
denote $\rho(\omega)=\tau_\omega(I_{\A(\omega)})$. Then
$\tau_\omega=\rho(\omega)\bar\tau_\omega$ where $\bar\tau_\omega$ is the
normalizedtrace on
the factor $\A(\omega)$. Clearly $\rho\in L^1(\omega,\mu)$, $\rho\geq 0$.
The formula (5.6)
can be rewritten  as follows:

$$\tau(a)=\int_\Omega \rho(\omega)\bar\tau_\omega(a(\omega))\mu(d\omega),
\ \ a\in\A\;.
\tag 5.7$$
The trace $\tau$
is normalized if and only if $\int_\Omega \rho(\omega)\mu(d\omega)=1$.

\proclaim {Lemma 5.3}
$M$ is finitely generated if and only if
$$\operatorname{ess\,sup}_{\omega\in\Omega}
(\rho(\omega)^{-1}\dim_{\tau(\omega)}M(\omega))<\infty\;.
\tag 5.8$$
\endproclaim

\subheading{Proof}
We refer to Theorem 11 in Chapter III of [Sc] for the following fact.
Assume that
$$M=\int_\Omega M(\omega)\mu(d\omega),\
M_1=\int_\Omega M_1(\omega)\mu(d\omega), \tag 5.9$$
are two Hilbert $\A$-modules decomposed into direct integrals according to
the decomposition (5.1) of the algebra $\A$. A continuous linear imbedding of
$\A$-modules $M\subset M_1$ exists iff
$$\dim_{\tau(\omega)}M(\omega)\leq\dim_{\tau(\omega)}M_1(\omega)
\tag 5.10$$
for almost all $\omega$. Applying this to compare the decomposition of
$M$ with the decomposition of $l^2(\A)\otimes\C^N$ we obtain the statement of
the Lemma.
$\square$

\medskip
{\bf Proof of (iv) in the general case.}
 Note first that it follows from the decomposition of the commutant (5.4) that
any representation of a discrete countable group
$\Ga$ in $M$  by automorphisms of $M$ as a Hilbert $\A$-module can be decomposed
into a direct integral of actions of $\Ga$ in $M(\om)$. This leads to
a decomposition of any flat Hilbert bundle with  fiber $M$
over a compact manifold $X$ into a direct integral of flat Hilbert bundles
with fibers $M(\om)$. This leads to the decomposition of the space
of $L^2$-forms with values in the Hilbert bundle, and furthermore to the
decomposition of the Laplacian $\Delta^{(p)}$ into a direct integral
of Laplacians $\Delta^{(p)}(\om)$ operating on forms with values in
flat vector bundles $E(\om)$ with  fibers $M(\om)$. The cocycle defining
the  bundle $E(\om)$ is obtained from   the components $g_{UV}(\om)$ in the
direct
decomposition of the cocycle $g_{UV}$ defining the bundle $E$.

All the functions of
the Laplacian $\De^{(p)}$ are then also decomposed. In particular, for the
spectral projection $E_\la^{(p)}$ we have the decomposition
$$E_\la^{(p)}=\int_\Om E_\la^{(p)}(\om)d\mu(\om)\;,\tag 5.11$$
where $ E_\la^{(p)}(\om)$ is the spectral projection of
$\De^{(p)}(\om)$. Taking into account Lemma 5.3 we see that the statement
(iv) of  Theorem 5.1 will be proved if we establish the estimate
$$\operatorname{ess\,sup}_{\om\in\Om}(\rho(\om)^{-1}\tau_\om(E_\la^{(p)}
(\om)))<\infty \tag 5.12$$
for any $\la>0$.

To prove this note first that for the Schwartz kernel
$E_\la^{(p)}(\om)(x,y)$ of the spectral projection
$E_\la^{(j)}(\om)$
we have an estimate
$$\text{ess sup}_{\om\in\Om}\sup_{x,y\in X}
\Vert E_\la^{(j)}(\om)(x,y)\Vert\leq C<\infty\;.
\tag 5.13$$
This follows from standard elliptic estimates (see e.g. arguments given in
\cite{Sh, FS, CG},
for similar estimates). The constant  $C$ in (5.13) does not depend on
$\om$ (but may depend on $\la$ and on all the other data, e.g. $X, \A$ and $E$).

Indeed the $L^2$-operator
norm of the orthogonal projection equals 1. Using the  elliptic  estimates
we can now prove  that
in fact the norm of the projection is bounded not only in $L^2$ but also
in the Sobolev spaces.
More precisely the  norm of $E^{(p)}_\la$ as an operator
from $H^{-s}$ to $H^s$ for any $s>0$ is bounded,
 so that the operator is infinitely smoothing.
 The constants in
these estimates depend on bounds for the coefficients, their derivatives
and  the
inverted symbol (the last is to insure the uniform ellipticity).
Therefore these estimates are uniform with respect to $\om$.
 Using them together with
the Sobolev embedding theorem ($H^s\subset C$ for $s>\dim X/2$ which
implies that the
Dirac delta-function belongs to $H^{-s}$ for $s>\dim X/2$) we arrive to (5.13).

Now note that due to Lemma 5.3
$$\operatorname{ess\,sup}_{\om\in\Om}(\rho(\om)^{-1}\tau_\om(I_\om))=
\operatorname{ess\,sup}_{\omega\in\Omega}
(\rho(\omega)^{-1}\dim_{\tau_\omega}M(\omega))<\infty$$
 because the fiber $M$ is finitely generated.
(Here $I_\om$ is the identity  operator in $M(\om)$.)

Using the estimate
$$|\tau_\omega(T)|\leq \Vert T\Vert\tau_\omega(I_\om),$$
which is true for any  endomorphism $T$ of the Hilbert module $M(\omega)$,
we immediately arrive at the estimate (5.2). This ends the proof of (iv)
in Theorem 5.1.
$\square$

\medskip
{\bf Remark.} Let $\A$ be a finite factor and $M$  a Hilbert $\A$-module. Then
$M$ is finitely generated if and only if $\dim_\tau M<\infty$. But this is
not true
for general finite von Neumann algebras. The following example is due to
M.Farber.
Consider  $\A=L^\infty([0,1])$
(where $[0,1]$ is taken with the Lebesgue measure) with the trace $\tau$ given
by the Lebesgue integral, and take
$$M=L^2([0,1/2])\oplus(L^2([1/2,3/4])\otimes \C^2)
\oplus(L^2([3/4,7/8])\otimes \C^3)\oplus\dots\,.$$
Here $a\in L^\infty([0,1])$ acts in each space
$L^2([(2^{k-1}-1)/2^{k-1},(2^k-1)/2^k]\otimes \C^k$ as
$a_k\otimes I_k$ where $a_k$ is the restriction of $a$ to
$[(2^{k-1}-1)/2^{k-1},(2^k-1)/2^k]$ and $I_k$ is the identity operator in
$\C^k$.
Clearly
$$\dim_\tau M=\sum_{k=1}^\infty{k\over 2^k}<\infty,$$
but $M$ is not finitely generated.


\Refs
\widestnumber\key {BFKM}

\ref\key A\by M.F.Atiyah\paper Elliptic operators, discrete groups and von
Neumann algebras\jour Ast\'erisque\vol 32-33\yr 1976\pages 43-72\endref

\ref\key BZ \by J.-M.Bismut, W.Zhang \paper An extension of a theorem by
Cheeger and M\"uller \jour Ast\'erisque \vol 205 \yr 1992 \endref

\ref\key BZ1 \by J.-M.Bismut, W.Zhang\paper Milnor and Ray-Singer metrics on
the equivariant determinant of a flat vector bundle
\jour Geometric and Functional Analysis \vol 4 \yr 1994 \pages 136-212 \endref

\ref\key BFKM \by D.Burghelea, L.Friedlander, T.Kappeler, P.McDonald
\paper Analytic and Reidemeister torsion for representations in finite type
Hilbert modules \jour Preprint \yr 1994
\endref

\ref \key CG\by J.Cheeger, M.Gromov\paper Bounds on the von Neumann dimension
of $L_2$-cohomology and the Gauss-Bonnet theorem for open manifolds
\jour J. Differential Geometry\vol 21\yr 1985\pages 1-34
\endref


\ref\key Di\by J.Dixmier\book Von Neumann algebras\yr 1981
\publ North-Holland Publishing Company
\endref

\ref \key D\by J. Dodziuk\paper De Rham-Hodge theory for $L^2$-cohomology
of infinite coverings\jour Topology\vol 16\pages 157-165\yr 1977
\endref

\ref \key E\by A.Efremov\paper Combinatorial and analytic Novikov-Shubin
invariants\jour Preprint
\endref

\ref \key E1\by A.Efremov\paper Cell decompositions and the Novikov-Shubin
invariants\jour Russ. Math. Surveys \vol 46\pages 219-220\yr 1991
\endref

\ref \key Fa\by M.Farber\paper Homological algebra of Novikov-Shubin
invariants and Morse inequalities\jour Prep\-rint\yr 1995\endref

\ref \key Fa1\by M.Farber\paper Von Neumann categories and extended
$L^2$ cohomology \jour Prep\-rint\yr 1996\endref

\ref \key FS \by B.V.Fedosov, M.A.Shubin \paper Index of random elliptic
operators I \jour Math. USSR Sbornik \vol 34 \yr 1978 \pages 671-699
\endref

\ref\key Fr\by P.Freyd\book Representations in Abelian Categories
\publ Proceedings of the Conference on Categorical Algebra, La Jolla 1965,
editors  S.Eilenberg, D.K.Harrison, S.MacLane, H.R\"ohrl,
Springer-Verlag \pages 95-120 \yr 1966
\endref


\ref\key GS\by M.Gromov, M.A.Shubin\paper Von Neumann spectra near
zero\jour Geometric and Functional Analysis\vol 1\yr 1991\pages 375-404\endref

\ref \key GS1\by M.Gromov, M.A.Shubin\paper Near-cohomology of Hilbert
complexes and topology of non-simp\-ly connected manifolds\jour Ast\'erisque
\vol 210\pages 283-294\yr 1992
\endref

\ref \key LL\by J.Lott, W.L\"uck\paper $L^2$-topological invariants of
3-manifolds\jour Invent. math. \yr 1995\vol 120\pages
\newline
15-60
\endref

\ref\key L\"u \by W.L\"uck \paper Hilbert modules and modules over finite
von Neumann algebras and applications to $L^2$-invariants
\jour Preprint Johannes Gutenberg-Universit\"at Mainz \yr 1995
\endref



\ref\key NS\by S.P.Novikov, M.A.Shubin\paper Morse inequalities and
von Neumann $II_1$-factors\jour Soviet
\newline
Math. Dokl. \vol 34:1 \yr 1987
\pages 79-82\endref

\ref \key NS1\by S.P.Novikov, M.A.Shubin\paper Morse inequalities and
von Neumann invariants of non\-simp\-ly connected manifolds\jour Uspehi
Matem. Nauk \vol 41:5\yr 1986\pages 222-223\endref

\ref \key P \by G.K.Pedersen \book $C^\ast$-algebras and their automorphism
groups \publ Academic Press \yr 1979 \endref

\ref \key Sc \by J.T.Schwartz \book $W^\ast$-algebras \publ Gordon and  Breach
\yr 1967 \endref

\ref\key Sh \by M.A.Shubin \paper Pseudodifferential almost-periodic
operators and von Neumann algebras
\newline
\jour Trans. Moscow Math. Soc.,
\yr 1976 \vol 35 \pages 103-163
\endref

\ref\key Sh1 \by M.A.Shubin \paper Spectral theory of elliptic operators on
non-compact manifolds \jour Ast\'erisque \vol 207 \pages 35-108
\endref

\ref\key Sh2 \by M.A.Shubin\paper Semiclassical asymptotics on covering
manifolds and Morse inequalities\jour Geometric and Functional Analysis
\yr 1996 \vol 6 \pages 370-409
\endref

\ref\key T \by M.Takesaki \book Theory of operator algebras I
\publ Springer-Verlag \yr 1979
\endref

\endRefs
\enddocument